# Properties of a renewal process approximation for a spin market model

Muffasir Badshah[1], Dr. Robert Boyer[1], Dr. Theodore Theodosopoulos[2]

[1]Department of Mathematics, Drexel University, Philadelphia, PA
[2]Department of Decision Sciences and Department of Mathematics, Drexel University, Philadelphia, PA

## Abstract

In this short note we investigate the nature of the phase transitions in a spin market model as a function of the interaction strength between local and global effects. We find that the stochastic dynamics of this stylized market model exhibit a periodicity whose dependence on the coupling constant in the Ising-like Hamiltonian is robust to changes in the temperature and the size of the market.

## 1. Model Overview

Over the past four years, a lattice-based spin model for market microstructure, motivated by studies of solid-state physics (e.g. the Ising potential), was proposed by Bornholdt and collaborators [2,5]. This model is governed by competition between interactions among agents at a local and global level. The resulting frustration generally leads to intermittent transitions between two phases: an "ordered" one, characterized by quasi-stable spin configurations, and a "disordered" one that is dominated by large-scale spin fluctuations.

Recent work by our group [1] explored the effects of variable neighborhood topology on the interaction Hamiltonian of the Bornholdt spin model. In the present paper we extend the phase transition analysis in [1] to construct a renewal process approximation for this spin market model. We show that the resulting renewal period is minimized for a value of the interaction strength $\alpha$ that is invariant with respect to changes in the inverse temperature $\beta$ and the number of interacting traders $N$.

The Bornholdt model provides a schematic environment to investigate high-frequency market. The investment attitude of the interacting trader is represented by the random variable, $S_i = \pm 1$, where +1 and –1 are indicative of a buyer and seller during a given period, respectively.

$$S_i(t+1) = +1 \quad \text{with} \quad p = \frac{1}{1+e^{-2\beta h_i(t)}}$$

$$S_i(t+1) = -1 \quad \text{with} \quad 1-p$$

$$h_i(t) = \sum_{j=1}^{N} J_{ij} S_j - \alpha S_i(t) |M(t)|$$

$$\text{where } M(t) = \frac{1}{N} \sum_{i=1}^{N} S_i(t)$$

$$\text{and } J_{ij} = \begin{cases} 1 & \text{if } j \text{ is a neighbor of } i \\ 0 & \text{otherwise} \end{cases}$$

The first term of the Hamiltonian captures the nearest neighbor interactions (effects of variations in the local neighborhood have been analyzed in [1]), with $J_{ij}$ modeling the neighborhood structure.

Agents in this model balance the desire to conform to the local majority, with the advantage of being in the global minority. In this sense, this model is a hybrid between the voter model [4] and the minority game [3]. It is plausible to interpret the Hamiltonian as a market driver in a variety of regimes. For instance, the third author showed that, asymptotically, a sufficiently strong coupling constant ($\alpha$) allows persistent market activity in a closely related spin market model [6], even in the frozen regime, i.e. when $\beta \to \infty$. Nevertheless, the characterization of the market behavior for intermediate values of these parameters remains to be explored. This is a central consideration for this paper.

## 2. Simulation Results

We simulated the dynamics of the Bornholdt model using a range of values for the coupling constant ($\alpha = 2, 3, \ldots, 8$), at decreasing temperatures ($T = \beta^{-1} = 5, 4, 3, 2$) and two different market sizes ($N = 16$ and $N = 32$).

We extend here the methodology we introduced in [1] for characterizing the paths of the stochastic market dynamics by studying the transition epochs between the "ordered" and "disordered" phases [5]. Let $\pi_{ord}$ denote the asymptotic percentage of time the system spends in the "ordered" phase. In [1] we showed that $\pi_{ord}$ is an increasing function of the connectivity in the local neighborhood. The local neighborhood used throughout the current study has degree eight. It can be seen in Figure 1 that there exists a unique value $\alpha^*$ of the interaction strength that minimizes $\pi_{ord}$. Indeed, this property of $\pi_{ord}$ remains true for all cases we investigated, irrespective of the temperature or market size.

Furthermore, let $\lambda_{ord}$ denote the Poisson rate that characterizes the transitions from the "ordered" to the "disordered" phase. In [1] we showed that $\lambda_{ord}$ is also an increasing function of the local connectivity. Figure 2 below shows that there exists a unique value $\hat{\alpha}$ of the interaction strength that maximizes $\lambda_{ord}$. Indeed, this property of $\lambda_{ord}$ remains true for all cases we investigated, irrespective of the temperature or market size.

Following the discussion in [1], we can approximate the paths of the local Hamiltonian in this market model as a renewal process, transitioning between the "ordered" and "disordered" phases. The expected duration of each sojourn in the "ordered" phase will last $\lambda_{ord}^{-1}$ time steps. On the other hand, the steady-state equation for this two-state system implies that

$$\pi_{ord} = \frac{\lambda_{dis}}{\lambda_{ord} + \lambda_{dis}}$$

and therefore the expected duration of each sojourn in the "disordered" phase will last $\lambda_{ord}^{-1}(\pi_{ord}^{-1} - 1)$ time steps. Thus, the renewal period $T_{renew}$ which characterizes the stochastic periodicity of this spin market model is given by $\frac{1}{\lambda_{ord}\pi_{ord}}$.

Figure 3 shows the dependence of the resulting renewal period on the coupling constant for different market sizes at different temperatures. Qualitatively we observe that all the curves in Figure 3 attain their unique minimum around the same value of $\alpha$. In order to test this hypothesis, we constructed spline models of $T_{renew}(\alpha)$ for each combination of $\beta$ and $N$. These models are shown in Figures 4a and 4b in the case of $N = 16$ and $N = 32$ respectively. Each of these models is minimized, and the resulting minima are shown as squares in Figures 3, 4a and 4b. As we can see, our simulation results allow us to conclude that $\alpha^* \cong \hat{\alpha} \cong 6.12 \pm 0.27 \ (p = 0.01)$.

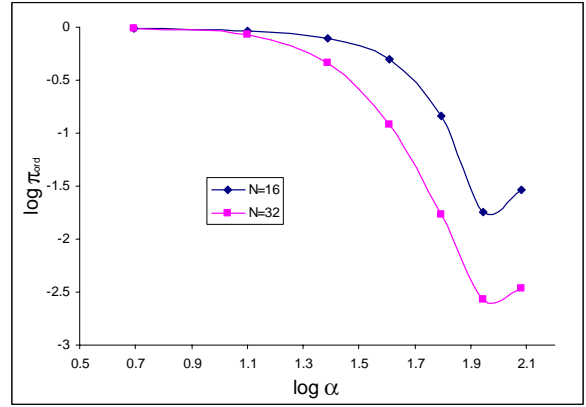

Figure 1: This graph shows $\pi_{ord}$ as a function of $\alpha$ for the case $\beta = 0.5$. Graphs for other temperatures are qualitatively similar (not included).

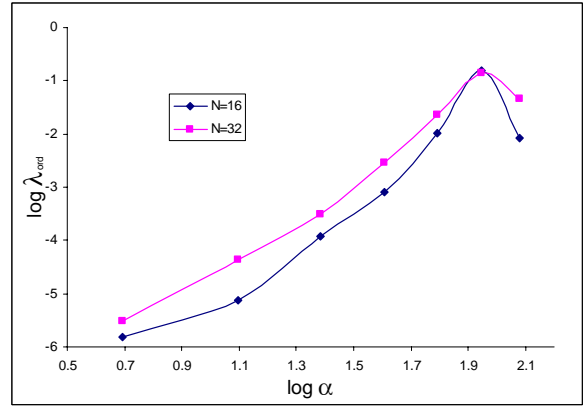

Figure 2: This graph shows $\lambda_{ord}$ as a function of $\alpha$ for the case $\beta = 0.5$. Graphs for other temperatures are qualitatively similar (not included).

Thus we conclude that there exists a critical value of the coupling constant which minimizes the stability of the spin market model, in the sense that it leads to the most frequent renewals. We proceeded to investigate the dependence of the resulting minimum renewal period on temperature and market size. Figure 5 shows that $T_{renew}(\alpha^*) \propto \beta^{\frac{2}{3}}$. In fact, the proportionality constant for the $N = 32$ case is about 1.5 times higher than the one for $N = 16$.

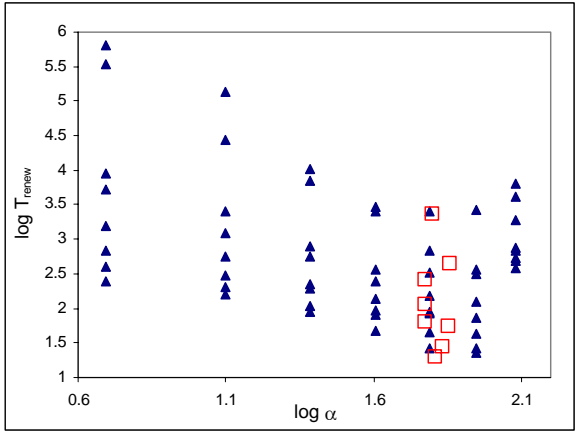

Figure 3: This graph shows $T_{renew}$ as a function of $\alpha$ for $T = \beta^{-1} = 5,4,3,2$ including both $N=16$ and $N=32$ cases.

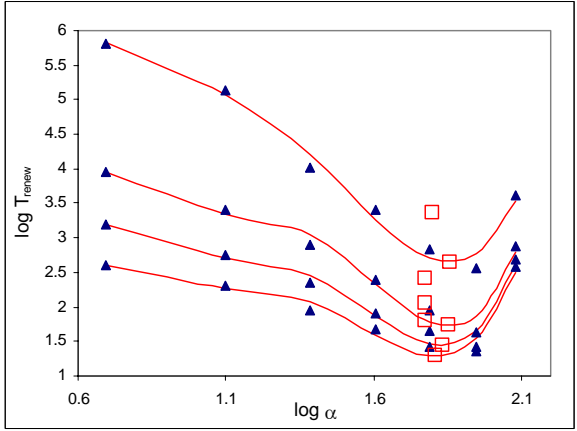

Figure 4a: Spline models for $T_{renew}(\alpha)$ for $N=16$.

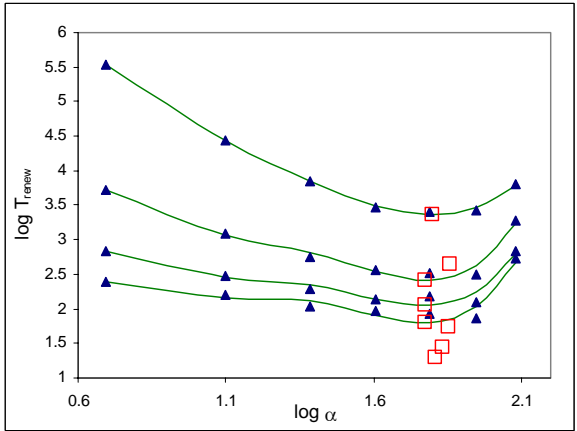

Figure 4b: Spline models for $T_{renew}(\alpha)$ for $N=32$.

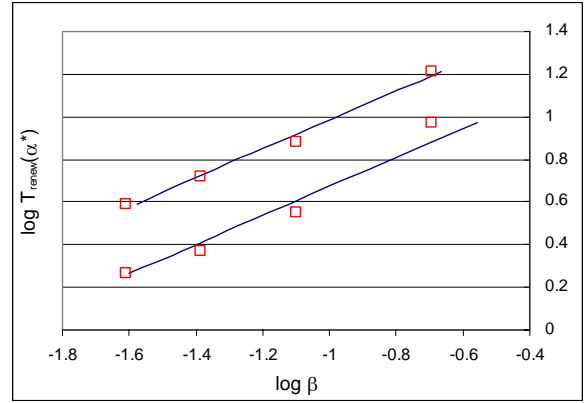

Figure 5: This graph shows $T_{renew}(\alpha^*)$ as a function of $\beta$ when $N=16$ (lower line) and $N=32$ (upper line).

## 3. Conclusion and Future Work

In summary we constructed a renewal process approximation to the stochastic dynamics of Bornholdt's spin market model and investigated the dependence of the resulting renewal period on the coupling constant at different intermediate temperatures and market sizes. Our simulation results lead us to formulate the following conjectures about this stylized market system:

- There exists a unique critical value of the coupling constant, approximately equal to 6.12, which maximizes the rate of transition from the "ordered" to the "disordered" phase and minimizes the asymptotic percentage of time spent in the "ordered" phase as well as the renewal period of the dynamics.
- This critical value is a constant characterizing the dynamics of the spin market model, and it does not depend on the temperature or the number of interacting agents.
- The minimum renewal period for this spin market model is proportional to the two thirds power of the inverse temperature, while the proportionality constant is an increasing function of the market size.

The critical value of the interaction strength we arrived at in this paper appears related to the value of the coupling constant that was shown in [6] to lead to the supercritical regime of a closely related process to the one we investigate here. Specifically in [6] the Bornholdt model was imbedded in an interacting particle system and a version of rapid stirring was applied to arrive at an analytic characterization of the resulting asymptotic dynamics. The investigation of any fundamental relationships between the critical value of $\alpha$ arrived at in the current paper and that formulated in [6] remains an open problem.

Furthermore, it is natural to extend the questions investigated here to the random neighborhood structures we presented in [1]. It remains plausible that $\alpha^*$ depends on the network topology, perhaps through the degree distribution, as was shown to be the case for $\pi_{ord}$ and $\lambda_{ord}$ in [1].

Finally, the nature of the sign change in $\frac{\partial T_{renew}}{\partial \alpha}$ at the critical value of the coupling constant remains to be understood as part of the market dynamics. It appears plausible that this is a qualitatively new effect, arising from the hybridization of two models (the voter model and the minority game) that do not possess it.